\documentclass[a4paper,12pt]{article}
\usepackage{latexsym}
\usepackage{amsmath}
\usepackage{amssymb}
\usepackage{amscd}
\usepackage{array}
\usepackage[all]{xy}
\usepackage{amsthm}
\usepackage{amsfonts}
\usepackage{enumerate}
\usepackage{multirow}

\usepackage{pdflscape}
\newtheorem{theorem}{Theorem}[]
\newtheorem{proposition}[theorem]{Proposition}

\newtheorem{corollary}[theorem]{Corollary}

\theoremstyle{definition}

\newtheorem{remark}[theorem]{Remark}

\newcommand{\Z}{\mathbf Z}

\newcommand{\F}{\mathbf F}

\newcommand{\Gal}{\mathrm{Gal}}

\newcommand{\Hol}{\mathrm{Hol}}

\newcommand{\Sym}{\operatorname{Sym}}

\newcommand{\GL}{\mathrm{GL}}

\newcommand{\End}{\operatorname{End}}

\newcommand{\Aut}{\operatorname{Aut}}

\newcommand{\wL} {{\widetilde{L}}}

\begin{document}
\begin{center}

\Large
Computation of Hopf Galois structures on low degree separable extensions and classification of those for degrees $p^2$ and $2p$

\vspace{1cm}
\large
Teresa Crespo and Marta Salguero

\vspace{0.3cm}
\footnotesize

Departament de Matem\`atiques i Inform\`atica, Universitat de Barcelona (UB), Gran Via de les
Corts Catalanes 585, E-08007 Barcelona, Spain, e-mail: teresa.crespo@ub.edu, msalguga11@alumnes.ub.edu

\end{center}

\date{\today}

\let\thefootnote\relax\footnotetext{Both authors acknowledge support by grant MTM2015-66716-P (MINECO/FEDER, UE).}

\normalsize
\begin{abstract} A Hopf Galois structure on a finite field extension $L/K$ is a pair $(H,\mu)$, where $H$ is a finite cocommutative $K$-Hopf algebra and $\mu$ a Hopf action. In this paper we present a program written in the computational algebra system Magma which gives all Hopf Galois structures on separable field extensions of degree up to eleven and several properties of those. Besides, we exhibit several results on Hopf Galois structures inspired by the program output. We prove that if $(H,\mu)$ is an almost classically Hopf Galois structure, then it is the unique Hopf Galois structure with underlying Hopf algebra $H$, up to isomorphism. For $p$ an odd prime, we prove that a separable extension of degree $p^2$ may have only one type of Hopf Galois structure and determine those of cyclic type; we determine as well the Hopf Galois structures on separable extensions of degree $2p$. We highlight the richness of the results obtained for extensions of degree 8 by computing an explicit example and presenting some tables which summarizes these results.

\end{abstract}

\section{Introduction}
A Hopf Galois structure on a finite extension of fields $L/K$ is a pair $(H,\mu)$, where $H$ is
a finite cocommutative $K$-Hopf algebra  and $\mu$ is a
Hopf action of $H$ on $L$, i.e a $K$-linear map $\mu: H \to
\End_K(L)$ giving $L$ a left $H$-module algebra structure and inducing a $K$-vector space isomorphism $L\otimes_K H\to\End_K(L)$.
Hopf Galois structures were introduced by Chase and Sweedler in \cite{C-S}.
For separable field extensions, Greither and
Pareigis \cite{G-P} give the following group-theoretic
equivalent condition to the existence of a Hopf Galois structure.

\begin{theorem}\label{G-P}
Let $L/K$ be a separable field extension of degree $g$, $\wL$ its Galois closure, $G=\Gal(\wL/K), G'=\Gal(\wL/L)$. Then there is a bijective correspondence
between the set of Hopf Galois structures on $L/K$ and the set of
regular subgroups $N$ of the symmetric group $S_g$ normalized by $\lambda (G)$, where
$\lambda:G \hookrightarrow S_g$ is the monomorphism given by the action of
$G$ on the left cosets $G/G'$.
\end{theorem}

For a given Hopf Galois structure on a separable field extension $L/K$ of degree $g$, we will refer to the isomorphism class of the corresponding group $N$ as the type of the Hopf Galois
structure. The Hopf algebra $H$ corresponding to a regular subgroup $N$ of $S_g$ normalized by $\lambda (G)$ is the sub-$K$-Hopf algebra $\wL[N]^G$ of the group algebra $\wL[N]$ fixed under the action of $G$, where $G$ acts on $\wL$ by $K$-automorphisms and on $N$ by conjugation through $\lambda$. The Hopf action is induced by $n \mapsto n^{-1}(\overline{1})$, for $n \in N$, where we identify $S_g$ with the group of permutations of $G/G'$ and $\overline{1}$ denotes the class of $1_G$ in $G/G'$. It is known that the sub-Hopf algebras of $\wL[N]^G$ are in 1-to-1 correspondence with the subgroups of $N$ stable under the action of $G$ (see e.g. \cite{CRV} Proposition 2.2) and that, given two regular subgroups $N_1, N_2$ of $S_g$ normalized by $\lambda (G)$, the Hopf algebras $\wL[N_1]^G$ and $\wL[N_2]^G$ are isomorphic if and only if the groups $N_1$ and $N_2$ are $G$-isomorphic.

Childs \cite{Ch1} gives an equivalent  condition to the existence of a Hopf Galois structure introducing the holomorph of the regular subgroup $N$ of $S_g$. We state the more precise formulation of this result due to Byott \cite{B} (see also \cite{Ch2} Theorem 7.3).

\begin{theorem}\label{theoB} Let $G$ be a finite group, $G'\subset G$ a subgroup and $\lambda:G\to \Sym(G/G')$ the morphism given by the action of
$G$ on the left cosets $G/G'$.
Let $N$ be a group of
order $[G:G']$ with identity element $e_N$. Then there is a
bijection between
$$
{\cal N}=\{\alpha:N\hookrightarrow \Sym(G/G') \mbox{ such that
}\alpha (N)\mbox{ is regular}\}
$$
and
$$
{\cal G}=\{\beta:G\hookrightarrow \Sym(N) \mbox{ such that }\beta
(G')\mbox{ is the stabilizer of } e_N\}
$$
Under this bijection, if $\alpha\in {\cal N}$ corresponds to
$\beta\in {\cal G}$, then $\alpha(N)$ is normalized by
$\lambda(G)$ if and only if $\beta(G)$ is contained in the
holomorph $\Hol(N)$ of $N$.
\end{theorem}

As a corollary to the preceding theorem Byott \cite{B}, Proposition 1, obtains the following formula to count Hopf Galois structures.

\begin{corollary}\label{cor} Let $L/K$ be a separable field extension of degree $g$, $\wL$ its Galois closure, $G=\Gal(\wL/K), G'=\Gal(\wL/L)$. Let $N$ be an abstract group of order $g$ and let $\Hol(N)$ denote the holomorph of $N$. The number $a(N,L/K)$ of Hopf Galois structures of type $N$ on $L/K$ is given by the following formula

$$a(N,L/K)= \dfrac{|\Aut(G,G')|}{|\Aut(N)|} \, b(N,G,G')$$

\noindent where $\Aut(G,G')$ denotes the group of automorphisms of $G$ taking $G'$ to $G'$, $\Aut(N)$ denotes the group of automorphisms of $N$ and $b(N,G,G')$ denotes the number of subgroups $G^*$ of $\Hol(N)$ such that there is an isomorphism from $G$ to $G^*$ taking $G'$ to the stabilizer in $G^*$ of $1_N$.

\end{corollary}

In Hopf Galois theory one has the following Galois correspondence theorem.

\begin{theorem}[\cite{C-S} Theorem 7.6]\label{esto} Let $(H,\mu)$ be a Hopf Galois structure on the field extension $L/K$.
For a sub-$K$-Hopf algebra $H'$ of $H$ we define
$$
L^{H'}=\{x\in L \mid \mu(h)(x)=\varepsilon(h)\cdot x \mbox{ for all } h\in H'\},
$$
where $\varepsilon$ is the counit of $H$.
Then, $L^{H'}$ is a subfield of $L$, containing $K$, and
$$
\begin{array}{rcl}
{\mathcal F}_H:\{H'\subseteq H \mbox{ sub-Hopf algebra}\}&\longrightarrow&\{\mbox{Fields }E\mid K\subseteq E\subseteq L\}\\
H'&\to &L^{H'}
\end{array}
$$
is injective and inclusion reversing.
\end{theorem}

In \cite{G-P} a class of Hopf Galois structures is identified for which the Galois correspondence is bijective. We shall say that a Hopf Galois structure $(H,\mu)$ on
$L/K$ is an \emph{almost classically Galois structure} if the corresponding regular subgroup $N$ of $S_g$ normalized by $\lambda(G)$ has the property that its centralizer $Z_{S_g}(N)$ in $S_g$ is contained in $\lambda(G)$.

\begin{theorem}[\cite{G-P} 5.2]
If $(H,\mu)$ is an almost classically Galois Hopf Galois structure on $L/K$, then the map ${\mathcal F}_H$ from the set of sub-$K$-Hopf algebras of $H$ into the set of subfields of $L$ containing $K$ is bijective.
\end{theorem}

In \cite{CRV2} the Hopf Galois character of separable field extensions of degree up to 7 and of some subextensions of its normal closure has been determined. In \cite{CRV} Theorem 3.4, a family of extensions is given with no almost classically Galois structure but with a Hopf Galois structure for which the Galois correspondence is bijective. In \cite{CRV4} a degree 8 non-normal separable extension having two non-isomorphic Hopf Galois structures with isomorphic underlying
Hopf algebras is presented.

\vspace{0.3cm}
In this paper we present a program written in the computational algebra system Magma which determines all Hopf Galois structures of a separable field extension of a given degree $g$ and their corresponding type. It is effective up to degree 11 and uses the Magma database of transitive groups which derives from the classification given in \cite{Bu}. Moreover our program distinguishes almost classically Galois structures and decides for the remaining ones if the Galois correspondence is bijective. Finally it classifies the Hopf Galois structures in Hopf algebra isomorphism classes. In the case of prime degree, we obtain the results already found in \cite{Ch1} theorem 2 and \cite{P} theorem 5.2, namely that if $L/K$ is a separable field extension of prime degree and $\wL$ its Galois closure, then $L/K$ has a Hopf Galois structure if and only if $\Gal(\widetilde{L}/K)$ is solvable and, in this case, the Hopf Galois structure is unique. We note that the case of degree 8 is especially interesting since there are 5 groups of order 8, up to isomorphism. We detail the results obtained in this case in Tables \ref{d8-ext}, \ref{d8-ext-cont} and \ref{d8-ext-Giso}. By performing an analysis of the outputs of our program, we have deduced several general behaviours. In Section \ref{almost} we prove that an almost classically Hopf Galois structure stands alone in its Hopf algebra isomorphism class. In Section \ref{p2} we prove that a separable field extension of degree $p^2$, for $p$ an odd prime, has at most one type of Hopf Galois structure and describe the ones of cyclic type. In Section \ref{2p} we determine the  Hopf Galois structures on separable field extensions of degree $2p$, for $p$ an odd prime.

\section{Description of the computation procedure}

Given a separable field extension $L/K$ of degree $g$, $\wL$ its Galois closure, $G=\Gal(\wL/K),$ \newline $G'=\Gal(\wL/L)$, the action of
$G$ on the left cosets $G/G'$ is transitive, hence the morphism $\lambda:G \rightarrow S_g$ identifies $G$ with a transitive subgroup of $S_g$, which is determined up to conjugacy. Moreover, if we enumerate the left cosets $G/G'$ starting with the one containing $1_G$, $\lambda(G')$ is equal to the stabilizer of $1$ in $G$. Therefore considering all separable field extensions $L/K$ of degree $g$ is equivalent to considering all transitive groups $G$ of degree $g$, up to conjugation.
The structure of the computation procedure is as follows:

\begin{enumerate}[{Step} 1]
\item Given a transitive group $G$ of degree $g$ and a type of regular subgroups $N$ of $S_g$, run over the conjugacy class of $N$ in $S_g$ and determine whether $N$ is normalized by $G$. In the affirmative case, check if the centralizer $Z(N)$ of $N$ in $S_g$ is contained in $G$. If it is so, the Hopf Galois structure determined by $N$ is almost classically Galois.
\item For each transitive group $G$ of degree $g$ and $G'=Stab_G(1)$, determine the number $intfields(G)$ of subgroups of $G$ containing $G'$, that is, by the fundamental  theorem of classical Galois theory, the number of intermediate fields of the extension $L/K$.
\item For each pair $(G,N)$ determined in Step 1, determine the number $subGst(N)$ of $G$-stable subgroups of $N$, i.e. subgroups of $N$ normalized by $G$, that is, the cardinality of the image of the map $\mathcal{F}_H$ in Theorem \ref{esto} for the Hopf Galois structure given by $N$. Check if this number equals $intfields(G)$, that is if the Galois correspondence is bijective.
\item For each pair $(G,N_1),(G,N_2)$, with $N_1\simeq N_2$ and  $subGst(N_1)=subGst(N_2)$, check if $N_1$ and $N_2$ are $G$-isomorphic, that is if the corresponding Hopf algebras are isomorphic. To this end, we use that for a regular subgroup $N$ of the symmetric group $S_g$, the automorphism group $\Aut(N)$ of $N$ is isomorphic to the stabilizer of 1 in the holomorph $Hol(N)$ of $N$ and that $Hol(N)$ is the normalizer of $N$ in $S_g$. We obtain the set of all isomorphisms by composing the isomorphism from $N_1$ to $N_2$ given by Magma with each automorphism of $N_2$. We run over this set of isomorphisms and check for each element whether it is a $G$-isomorphism until the answer is affirmative or the set is exhausted.
\end{enumerate}

We note that in Step 1 we compute the transversal of the normalizer of $N$ in $S_g$ and the conjugate of $N$ by each element in this transversal. This computation occurs to need a significantly shorter  execution time than the use of the Magma function Class from degree 9 onwards. The program returns all regular subgroups $N$ of $S_g$ giving a Hopf Galois structure, hence determines explicitly all of them. In the vector which collects such $N$'s we have added a numbering variable in order to identify each of them with an integer number. This numeration is respected all along the program so that, once the $N$'s have been computed in Step 1, we can easily know the properties of the corresponding Hopf Galois structures by searching the assigned number. This greatly simplifies the reading and interpretation of the results. The Magma code of this program may be found in \cite{CS}.

\section{Almost classically Galois Hopf Galois structures}\label{almost}

By looking at the distribution in  Hopf algebra isomorphism classes of the Hopf Galois structures of a given separable extension provided by our program we have deduced the following result.

\begin{proposition} Let $L/K$ be a separable field extension of degree $g$. Let $(\mathcal{H},\mu)$ be an almost classically Galois structure on $L/K$ and  $(\mathcal{H}',\mu')$ a Hopf Galois structure on $L/K$. If the Hopf algebras $\mathcal{H}$ and $\mathcal{H}'$ are $K$-isomorphic, then the  Hopf Galois structures $(\mathcal{H},\mu)$ and $(\mathcal{H}',\mu')$ are isomorphic. Hence an almost classically Galois structure stands alone in its Hopf algebra isomorphism class.
\end{proposition}

\noindent
{\it Proof.} Let $\wL$ be a Galois closure of $L/K$, $G=\Gal(\wL/K), G'=\Gal(\wL/L)$, let $\lambda:G \rightarrow S_g$ be the monomorphism from G into the symmetric group $S_g$ given by the action of $G$ on the left cosets $G/G'$. By Theorem \ref{G-P} and \cite{G-P} Proposition 4.1, $(\mathcal{H},\mu)$ corresponds to a regular subgroup $N$ of $S_g$, normalized by $\lambda(G)$ and such that the centralizer $Z$ of $N$ in $S_g$ is contained in $\lambda(G)$ and $(\mathcal{H}',\mu')$ corresponds to a regular subgroup $N'$ of $S_g$, normalized by $\lambda(G)$. We know also that an isomorphism from $\mathcal{H}$ to $\mathcal{H}'$ corresponds to a $\lambda(G)$-isomorphism from $N$ to $N'$. If $f$ is such an isomorphism, we have $\sigma f(n)\sigma^{-1}=f(\sigma n\sigma^{-1})$, for all $\sigma \in \lambda(G)$ and $n \in N$. Now, since $Z\subset \lambda(G)$, we have, for all $z \in Z, n \in N$, $zf(n)z^{-1}=f(znz^{-1})=f(n)$. This implies that $N'=f(N)$ is contained in the centralizer of $Z$ in $S_g$. Since $N$ is regular, this centralizer coincides with $N$. We have then $N'=N$ and, again by Theorem \ref{G-P}, this implies that the  Hopf Galois structures $(\mathcal{H},\mu)$ and $(\mathcal{H}',\mu')$ are isomorphic. $\Box$

\section{Extensions of degree $p^2$, for $p$ an odd prime}\label{p2}

For $p$ prime, there are exactly two groups of order $p^2$, up to isomorphism, the cyclic one $C_{p^2}$ and the direct product of two copies of $C_p$, hence two possible types for a Hopf Galois structure of a field extension of degree $p^2$. We shall prove that the two types do not occur simultaneously, when $p\neq 2$. This fact was suggested to us by the program output for degree 9 extensions. The case $p=2$ goes differently. Both Galois extensions of degree 4 and separable extensions of degree 4 whose Galois closure has Galois group the dihedral group $D_{2\cdot 4}$ have Hopf Galois structures of cyclic type and of type $C_2\times C_2$. If we write $C_{p^2}$ additively as $\Z/p^2 \Z$, its holomorph is $\Z/p^2 \Z \rtimes (\Z/p^2 \Z)^*$. For $C_p\times C_p$ the automorphism group is isomorphic to $\GL(2,\F_p)$.

\begin{proposition}\label{psqua}
Let $L/K$ be a separable field extension of degree $p^2$, $p$ an odd prime, $\widetilde{L}/K$ its normal closure and $G\simeq \Gal(\widetilde{L}/K)$.
If $L/K$ has a Hopf Galois structure of type $C_{p^2}$, then it has no structure of type $C_p\times C_p$. Therefore a separable field extension of degree $p^2$, $p$ an odd prime,  has at most one type of Hopf Galois structures either cyclic or elementary abelian.
\end{proposition}

\noindent {\it Proof.}
 By theorem \ref{theoB}, if $L/K$ has a Hopf Galois structure of type $C_{p^2}$, then $G$ is a transitive subgroup of $\Hol(C_{p^2})$. We shall see that all transitive subgroups of $\Hol(C_{p^2})$ contain an element of order $p^2$. Let us write $\Hol(C_{p^2})$ as $\Z/p^2 \Z \rtimes (\Z/p^2 \Z)^*$ and let $\sigma$ be a generator of $(\Z/p^2 \Z)^*$. The immersion of  $\Hol(C_{p^2})$ in the symmetric group $S_{p^2}$ is given by sending the generator $1$ of $\Z/p^2 \Z$ to the $p^2$-cycle $(1,2,\dots,p^2)$ and $\sigma$ to itself, considered as a permutation. The stabilizer of 0 in the image $H$ of  $\Hol(C_{p^2})$ in $S_{p^2}$ consists of the images of the elements $(0,\sigma^j)$. We have $|H|=|\Hol(C_{p^2})|=p^3(p-1)$, hence $H$ has a unique $p$-Sylow subgroup $Syl(H)$ which is isomorphic to the only non-abelian group of order $p^3$ having an element of order $p^2$ (see \cite{C}). Now, a subgroup $H'$ of $H$ is transitive if and only if $[H':Stab_H(0) \cap H']=p^2$. Let $H'$ be a transitive subgroup of $H$. We have then $p^2 \mid |H'|$ and $|H'| \mid p^3(p-1)$, hence $H'$ has a unique $p$-Sylow subgroup $Syl(H')$ which has order $p^3$ or $p^2$. In the first case, $Syl(H')=Syl(H)$ contains an element of order $p^2$. In the second case, $Syl(H')$ is a subgroup of $Syl(H)$ of order $p^2$. The group $Syl(H)$ is isomorphic to the group

 $$G_p:= \left\{ \left(\begin{array}{cc} 1+pm & b \\ 0 & 1 \end{array} \right) : m, b \in \Z/p^2 \Z \right\},$$

 \noindent where $m$ actually only matters modulo $p$. The group $G_p$ has $p^3-p^2$ elements of order $p^2$, those with $b\not \equiv 0 \pmod{p}$, hence $p$ cyclic subgroups of order $p^2$ and $p^2-1$ elements of order $p$, those nontrivial with $b \equiv 0 \pmod{p}$, hence one noncyclic subgroup of order $p^2$. Then $H'$ contains an element of order $p^2$ except in the case in which $Syl(H')$ is isomorphic to the noncyclic subgroup of order $p^2$ of $G_p$. The corresponding subgroup of $\Hol(C_{p^2})$ is generated by $(p,Id)$ and $(0,\sigma^{p-1})$. Its intersection with $Stab_H(0)$ consists in the elements $(0,\sigma^{\l(p-1)}), 1\leq \l \leq p$, hence this intersection has order $p$. We have then that if $Syl(H´)$ is isomorphic to the noncyclic subgroup of order $p^2$ of $G_p$, then $p$ divides exactly $[H':Stab_H(0) \cap H']$ and $H'$ is not transitive. We have proved then that all transitive subgroups of $\Hol(C_{p^2})$ contain an element of order $p^2$.

 Let us look now at $Hol(C_p\times C_p)$. By \cite{Ko}, Theorem 4.4, $Hol(C_p\times C_p)$ has no elements of order $p^2$. Taking into account what we have proved above, this finishes the proof of the proposition.
 $\Box$

\begin{remark} Kohl  proves in \cite{Ko} that any Hopf Galois structure on a cyclic extension of order $p^n$, for $p$ an odd prime, is of cyclic type. Childs studies in \cite{Ch3} these Hopf Galois structures in the case of cyclic extensions of order $p^2$.
\end{remark}

We give a more precise description of the Hopf Galois structures of cyclic type on separable field extensions of degree $p^2$ in the next theorem.

\begin{theorem}\label{p^2} Let $L/K$ be a separable field extension of degree $p^2$, $p$ an odd prime, $\widetilde{L}/K$ its normal closure and $G\simeq \Gal(\widetilde{L}/K)$.
 The extension $L/K$ has a Hopf Galois structure of cyclic type if and only if $G$ is isomorphic to the semidirect product $C_{p^2} \rtimes C_m$, for $m$ a divisor of $p(p-1)$. The number of structures is $p$ for $m=1$ and $m=p$ and is equal to $1$ in the remaining cases.
\end{theorem}

\noindent {\it Proof.} We have proved that all transitive subgroups of $\Hol(C_{p^2})$ have an element of order $p^2$. Reciprocally  a subgroup of $\Hol(C_{p^2})$ having an element of order $p^2$ is transitive. Let us write $\Hol(C_{p^2})=\langle 1, \sigma \rangle$, as above. The cyclic subgroups of order $p^2$ of $\Hol(C_{p^2})$ are $\langle (1,\sigma^{j(p-1)}) \rangle$, for $j=0,\dots, p-1$. We obtain then that the transitive subgroups of $\Hol(C_{p^2})$ are these $p$ groups of order $p^2$ and one group isomorphic to $C_{p^2} \rtimes C_m$, for each divisor $m$ of $p(p-1)$, $m\neq 1$, namely $\langle (1,\sigma^{j(p-1)}), (0,\sigma^{p(p-1)/m}) \rangle$. We count the number of structures using Corollary \ref{cor}. For $G\simeq C_{p^2}$, the number of structures is clearly equal to the number of transitive subgroups of $\Hol(C_{p^2})$ isomorphic to $C_{p^2}$, hence $p$. For $G=\langle a,b\rangle$, with $a=(1,\sigma^{j(p-1)}), b=(0,\sigma^{p(p-1)/m})$, we have $G'=\langle b \rangle$. An automorphism of $G$ sending $G'$ to $G'$ must be the identity on $G'$. Now the image of $a$ must be $a^i$, with $\gcd(i,p)=1$, when $m \neq p$, hence $|\Aut(G,G')|=|\Aut(N)|$ and the number of structures is 1. When $m=p$, the image of $a$ may be $a^ib^{(p-1)j}$, with $\gcd(i,p)=1$, and $0\leq j \leq p-1$. Hence $|\Aut(G,G')|=p^2(p-1)=p|\Aut(N)|$ and the number of structures is p. $\Box$

\section{Extensions of degree $2p$, for $p$ an odd prime}\label{2p}

In the next two theorems we determine the Hopf Galois structures of separable field extensions of degree $2p$, where $p$ is an odd prime. This result was suggested by the output of our program for degree 6 and 10 extensions. Since there are exactly two groups of order $2p$, up to isomorphism, the cyclic one and the dihedral one, we have exactly two types of Hopf Galois structures. We deal with the cyclic type in Theorem \ref{2p1} and with the dihedral type in Theorem \ref{2p2}. We note that the case of Galois extensions was already obtained by Byott in \cite{B2} and by Kohl in \cite{Ko2}.

\begin{theorem}\label{2p1} Let $L/K$ be a separable field extension of degree $2p$, $p$ an odd prime, $\widetilde{L}/K$ its normal closure and $G\simeq \Gal(\widetilde{L}/K)$.
 $L/K$ has a Hopf Galois structure of cyclic type if and only if $G$ is isomorphic to either

\begin{enumerate}[(i)]
\item the semidirect product $C_{2p} \rtimes C_m$ of a cyclic group of order $2p$ and a cyclic group of order $m$ dividing $p-1$
\item or the semidirect product $C_p \rtimes C_{m}$ of a cyclic group of order $p$ and a cyclic group of even order $m$ dividing $p-1$.
\end{enumerate}

The number of structures is 1, except in case (ii), for $m=1$, when this number is $p$. In this last case, $G$ is isomorphic to the dihedral group $D_{2p}$.
\end{theorem}

\noindent {\it Proof.}
By theorem \ref{theoB}, if $L/K$ has a Hopf Galois structure of type $C_{2p}$, then $G$ is a transitive subgroup of $\Hol(C_{2p})$. We have that

    $$\Hol(C_{2p}) \simeq C_{2p} \rtimes C_{p-1} = \langle a,b \mid a^{2p}=1, b^{p-1}=1, bab^{-1}=a^i \rangle,$$

\noindent where $i$ has order $p-1$ modulo $2p$. By ordering the elements in $C_{2p}$ as $a,a^2,\cdots,a^{2p}=1$, we obtain the embedding of $\Hol(C_{2p})$ in the symmetric group $S_{2p}$ mapping $a$ to $(1,2,\cdots,2p)$ and $b$ to a permutation of order $p-1$ sending any number in $\{1,2,\cdots,2p\}$ to one with the same parity. Let $G$ be a transitive group of degree $2p$. The order of $G$ is a multiple of $2p$, and $G$ contains an element of order $p$. Since the only subgroup of order $p$ of $\Hol(C_{2p})$ is $\langle a^2 \rangle$, then $G$ contains $a^2$.  Moreover, in order to be transitive, it must contain an element sending $2p$ to 1. Such elements in $\Hol(C_{2p})$ are exactly those of the form $a^jb^k$, with $j$ odd. Now $\langle a^2, a^jb^k \rangle=\langle a^2, ab^k \rangle$ if $j$ is odd. The transitive subgroups of $\Hol(C_{2p})$ containing $a$ are clearly $\langle a, b^{(p-1)/m} \rangle \simeq C_{2p} \rtimes C_m$, where $m$ ranges over the positive divisors of $p-1$. This gives case (i).

\vspace{0.3cm}
Let us assume now that $a \not \in G$. If $|G|=2pm$, then $G$ contains exactly $m$ elements sending $2p$ to 1. Let us determine the number of elements sending $2p$ to 1 in $G_k:=\langle a^2, ab^k \rangle$. Since $(ab^k)a^2(ab^k)^{-1}=a^{2i^k}$, we have $(ab^k)^l=a^{1+i^k+\cdots +i^{(l-1)k}} b^{kl}$ and we obtain that $G_k$ contains the elements $ab^{kl}$, with $l$ odd. If the order of $b^k$ is odd, then $a \in G_k$. If the order $n$ of $b^k$ is even, then $G_k$ contains $n/2$ elements sending $2p$ to 1, namely $ab^{kl}$, with $l$ odd, $1\leq j \leq n-1$. Moreover, the elements of order $n$ in $\langle b \rangle$ are among the elements $b^{kl}$, with $l$ odd, $1\leq j \leq n-1$. We have then that the transitive subgroups of $\Hol(C_{2p})$ not containing $a$ are  $\langle a^2, ab^{(p-1)/2m} \rangle \simeq C_{p} \rtimes C_{m}$, where $m$ ranges over the even positive divisors of $p-1$. This gives case (ii).

\vspace{0.3cm}
We determine now the number of structures by using Corollary \ref{cor}. For $N=C_{2p}$, we have $|\Aut(N)|=p-1$. We have exactly one subgroup of $\Hol(C_{2p})$ for each isomorphism class of $G$. For the groups $G=\langle a, b^{(p-1)/m}\rangle$, the stabilizer $G'$ of $2p$ is $\langle b^{(p-1)/m}\rangle$. The image of $a$ under an automorphism of $G$ is $a^j$, with $j$ coprime with $2p$. Since $bab^{-1}=a^i$, we have that an automorphism of $G$ sending $G'$ to $G'$ must send $b^{(p-1)/m}$ to itself. We obtain then $|\Aut(G,G')|=|\Aut(N)|$, hence the number of Hopf Galois structures is 1. For the groups $\langle a^2, ab^{(p-1)/2m} \rangle$,  the stabilizer $G'$ of $2p$ is $\langle b^{(p-1)/m}\rangle$. If $m=1$, $G \simeq D_{2p}$ and we obtain $|\Aut(G,G')|=p(p-1)$, hence the number of Hopf Galois structures is $p$. If $m>1$, taking into account the conjugate of $a^2$ by $ab^{(p-1)/2m}$, we obtain that an automorphism of $G$ sending $G'$ to $G'$ must send $ab^{(p-1)/2m}$ to itself, hence $|\Aut(G,G')|=|\Aut(N)|$ and the number of Hopf Galois structures is 1. $\Box$

\begin{theorem}\label{2p2} Let $L/K$ be a separable field extension of degree $2p$, $p$ an odd prime, $\widetilde{L}/K$ its normal closure and $G\simeq \Gal(\widetilde{L}/K)$. The extension $L/K$ has a Hopf Galois structure of dihedral type if and only if $G$ is isomorphic to either

\begin{enumerate}[(i)]
\item the semidirect product $(C_{p} \times C_p) \rtimes C_m$ of the direct product of two cyclic groups of order $p$ and a cyclic group of even order $m$ dividing $p-1$;
\item or the semidirect product $(C_p \times C_p)\rtimes (C_2 \times C_{m})$ of the direct product of two cyclic groups of order $p$ and the direct product of  a cyclic group of order 2 and a cyclic group of even order $m$ dividing $p-1$;
\item or the semidirect product $C_{2p}  \rtimes C_m$ of a cyclic group of order $2p$ and a cyclic group of order $m$ dividing $p-1$;
\item or the semidirect product $C_{p}  \rtimes C_m$ of a cyclic group of order $p$ and a cyclic group of even order $m$ dividing $p-1$.
\end{enumerate}

The number of structures is always 2.

\end{theorem}

\noindent {\it Proof.} By theorem \ref{theoB}, if $L/K$ has a Hopf Galois structure of type $D_{2p}$, then $G$ is a transitive subgroup of $\Hol(D_{2p})$. Let us write $D_{2p}=\langle \rho, \sigma \mid \rho^p=1, \sigma^2=1, \sigma \rho \sigma=\rho^{-1} \rangle$. We shall see that the automorphism group of $D_{2p}$ is isomorphic to $\Hol(C_p)$. More precisely, it is generated by $\varphi$ and $\psi$ determined by $\varphi(\rho)=\rho, \varphi(\sigma)=\sigma \rho$ and $\psi(\rho)=\rho^i, \psi(\sigma)=\sigma$, where $i$ is a generator of $(\Z/p\Z)^*$, and satisfying $\psi\varphi\psi^{-1}=\varphi^i$. We obtain then that $\Hol(D_{2p})$ is generated by $\rho, \sigma, \varphi, \psi$ with the relations

    $$ \begin{array}{ccccc} \rho^{p}=1, & \sigma^2=1, & \varphi^p=1, & \psi^{p-1}=1, & \sigma \rho \sigma=\rho^{-1}, \\ \psi\varphi\psi^{-1}=\varphi^i, & \sigma \varphi \sigma = \rho\varphi, & \psi\rho\psi^{-1}=\rho^i, & \rho\varphi=\varphi\rho, & \sigma \psi=\psi \sigma\end{array}$$

\noindent and is isomorphic to $(C_{p} \times C_p) \rtimes (C_2 \times C_{p-1})$.
By ordering the elements in $D_{2p}$ as $\rho,\rho^2,\cdots,\rho^{p}=1,\sigma\rho,\sigma\rho^2,\cdots,\sigma\rho^{p}=\sigma$, we obtain the embedding of $\Hol(D_{2p})$ in the symmetric group $S_{2p}$ mapping $\rho$ to $(1,2,\cdots,p)(2p,2p-1,\cdots,p+1)$, $\sigma$ to $(1,p+1)(2,p+2),\cdots (p,2p)$, $\varphi$ to $(p+1,p+2,\cdots,2p)$ and $\psi$ to $(1,i,i^2,\cdots,i^{p-2})(p+1,p+i,p+i^2,\cdots,p+i^{p-2})$, where the powers of $i$ are computed modulo $p$.

If $G$ is a transitive subgroup of $\Hol(D_{2p})$, then $2p$ divides the order $|G|$ of $G$, hence $G$ contains an element of order $p$ which belongs to the $p$-Sylow subgroup $\langle \rho,\varphi \rangle$ of $\Hol(D_{2p})$. We distinguish two cases, depending on wether the $p$-Sylow subgroup of $G$ has order $p$ or $p^2$.

\begin{enumerate}[]
\item {\bf Case 1.} If the $p$-Sylow subgroup of $G$ has order $p^2$, then it is $\langle \rho, \varphi \rangle$. If $G$ is transitive, it must contain a permutation sending $p$ to $2p$ and if $|G|=2mp$, then it contains $m$ elements sending $p$ to $2p$. The elements in $\Hol(D_{2p})$ sending $p$ to $2p$ are those of the form $\sigma x$, with $x \in Stab(p)=\langle \varphi, \psi \rangle$. We determine now the order of these elements. First the order of $\sigma$ is 2. If $0< j \leq p-1$, we have $(\sigma \varphi^j)^n=\rho^{jn/2}\varphi^{jn}$, if $n$ is even, and $(\sigma \varphi^j)^n=\sigma \rho^{j(n-1)/2}\varphi^{jn}$, if $n$ is odd. Hence the order of $\sigma \varphi^j$ is $2p$. If $0\leq j \leq p-1$ and $0<k<p-1$, we have $(\sigma \varphi^j \psi^k)^n=\rho^{j(\sum_{m=0}^{(n-2)/2}i^{2mk})}\varphi^{j(\sum_{m=0}^{n-1}i^{mk})}$, if $n$ is even, and $(\sigma \varphi^j \psi^k)^n=\sigma\rho^{j(\sum_{m=0}^{(n-3)/2}i^{(2m+1)k})}\varphi^{j(\sum_{m=0}^{n-1}i^{mk})}$, if $n$ is odd. Hence the order of $\sigma \varphi^j \psi^k$ is equal to the order of $\psi^k$ if this order is even and equal to twice the order of $\psi^k$ if this order is odd.

    We obtain then the following transitive subgroups of $\Hol(D_{2p})$ having a $p$-Sylow group of order $p^2$:

    \begin{enumerate}[a)]
    \item $\langle \varphi, \rho, \sigma, \psi^{2(p-1)/m} \rangle \simeq \langle \varphi, \rho, \sigma\psi^{(p-1)/2}, \psi^{2(p-1)/m} \rangle \simeq (C_p\times C_p)\rtimes C_m$, for $m$ a divisor of $p-1$, exactly divisible by 2.
     \item $\langle \varphi, \rho, \sigma \psi^{(p-1)/m} \rangle \simeq (C_p\times C_p)\rtimes C_m$, for $m$ a divisor of $p-1$, divisible by 4.
      \item $\langle \varphi, \rho, \sigma, \psi^{(p-1)/m} \rangle \simeq (C_p\times C_p)\rtimes (C_2 \times C_m)$, for $m$ an even divisor of $p-1$.
    \end{enumerate}

The groups of types a) and b) correspond to (i) in the statement and those of type c) correspond to (ii). We determine now the number of structures by using Corollary \ref{cor}. For $N=D_{2p}$, we have $|\Aut(N)|=p(p-1)$. The stabilizer $G'$ of $p$ is $\langle \varphi, \psi^k \rangle$, where $k=2(p-1)/m$ in cases a) and b) and $k=(p-1)/m$ in case c). We determine now the automorphisms of $G$ sending $G'$ to $G'$. Let $h$ be such an automorphism. The image of $\varphi$ under $h$ must be $\varphi^{j_1}$, with $1\leq j_1 \leq p-1$ and, taking into account the action of $\psi$ on $\varphi$, $h$ must send $\psi^k$ to $\varphi^{j_2} \psi^k$, with $0\leq j_2 \leq p-1$. The image of $\rho$ under $h$ must belong to $\langle \rho, \varphi \rangle$ and, since $\sigma$ has order 2 and commutes with $\psi$, the image of $\sigma$ under $h$ must be an order 2 element in $\langle \sigma, \psi \rangle$, hence $\sigma$, or $\sigma \psi^{(p-1)/2}$, in case $\psi^{(p-1)/2}$ belongs to $G$, i.e. for cases b) and c). Taking into account the action of $\sigma$ on $\varphi$ and $\rho$, we obtain that the elements in $\Aut(G,G')$ are of one of the following forms, with $1\leq j_1 \leq p-1, 0\leq j_2 \leq p-1$ and the second form not occurring in case a).

$$\begin{array}{lll} \varphi & \mapsto &\varphi^{j_1} \\
\psi & \mapsto & \varphi^{j_2} \psi  \\
\rho & \mapsto & \rho^{j_1}  \\
\sigma & \mapsto & \sigma  \end{array}
\quad \quad \quad
\begin{array}{llll} \varphi & \mapsto &\varphi^{j_1} \\
\psi & \mapsto & \varphi^{j_2} \psi  \\
\rho & \mapsto & \rho^{-j_1} \varphi^{-2j_1}  \\
\sigma & \mapsto & \sigma\psi^{(p-1)/2}  \end{array}
$$

We obtain then $|\Aut(G,G')|=|\Aut(N)|$ in case a) and $|\Aut(G,G')|=2|\Aut(N)|$ in cases b) and c). Hence in all three cases, the number of Hopf Galois structures is~2.

\item {\bf Case 2.} We consider now the case when the $p$-Sylow subgroup of $G$ has order $p$. The subgroups of order $p$ of $\Hol(D_{2p})$ are $\langle \rho\varphi^j \rangle$, with $0\leq j \leq p-1$, and $\langle \varphi \rangle$. To be transitive, $G$ must contain an element sending $p$ to $2p$, hence an element in $\sigma \langle \varphi, \psi \rangle$. Since the normalizers in $\Hol(D_{2p})$  of $\langle \rho\varphi^j \rangle$, if $j \not \in \{0,2\}$, and the one of $\langle \varphi \rangle$ reduce to $\langle \rho, \varphi, \psi \rangle$, the $p$-Sylow subgroup of $G$ is $\langle x \rangle$, where $x= \rho$ or $x=\rho\varphi^2$. If $G\simeq D_{2p}, G$ must contain an element $y$ of order 2 satisfying $yxy=x^{-1}$ and sending $p$ to $2p$. We obtain then the groups $\langle \rho, \sigma \rangle \simeq \langle \rho \varphi^2, \sigma \psi^{(p-1)/2}\rangle \simeq D_{2p}$. If $G\simeq C_{2p}$, $G$ must contain an element $y$ of order 2 sending $p$ to some element in $\{p+1,\cdots, 2p\}$ and commuting with $x$. We have $y=\sigma \psi^{(p-1)/2}(\rho \varphi^{2})^j$, for $x=\rho$ and $y=\sigma\rho^k$ for $x=\rho\varphi^2$. We obtain then the groups $\langle  \sigma \psi^{(p-1)/2} \rho^{j+1}\varphi^{2j} \rangle  \simeq\langle \sigma\rho^k \varphi  \rangle \simeq C_{2p}, 0\leq j \leq p-1, 0\leq k \leq p-1$. If a subgroup of $\Hol(D_{2p})$ contains one of these copies of $C_{2p}$, then it is transitive.

    By computation, we obtain that the only subgroups of $\Hol(D_{2p})$ strictly containing $\langle \sigma \psi^{(p-1)/2} \rho^{j+1}\varphi^{2j}  \rangle$ and having a $p$-Sylow subgroup of order $p$ are $\langle  \sigma \psi^{(p-1)/2} \rho^{j+1}\varphi^{2j} ,$ \linebreak
    $\varphi^{j(-i^l+1)/2}\psi^l  \rangle$, where $l=(p-1)/d$, for some divisor $d$ of $p-1$. Moreover $\langle  \sigma \psi^{(p-1)/2} \rho^{j+1}\varphi^{2j}  \rangle$ is normal in such a group and the order of $\langle  \sigma \psi^{(p-1)/2} \rho^{j+1}\varphi^{2j} ,
    \varphi^{j(-i^l+1)/2}\psi^l  \rangle$ is $2pd$. Analogously, the only subgroups of $\Hol(D_{2p})$ strictly containing $\langle \sigma\rho^k \varphi  \rangle$ and having a $p$-Sylow subgroup of order $p$ are $\langle \sigma\rho^k \varphi, \varphi^{(i^l-1)((p+1)/2-k)}\psi^l \rangle$, where $l=(p-1)/d$, for some divisor $d$ of $p-1$. Moreover $\langle \sigma\rho^k \varphi  \rangle$ is normal in such a group and the order of $\langle \sigma\rho^k \varphi, \varphi^{(i^l-1)((p+1)/2-k)}\psi^l \rangle$ is $2pd$. We have then $2p$ transitive subgroups of $\Hol(D_{2p})$ isomorphic to the semidirect product $C_{2p} \rtimes C_d$, for each divisor $d$ of $p-1$. This gives case (iii) in the statement.

    Finally, if $G$ has no element of order $2p$, it still must contain an element sending $p$ to $2p$. Such elements in $\Hol(D_{2p})$ are of the form $\sigma \varphi^k \psi^l$. The order of $\sigma \varphi^k \psi^l$ is equal to the order of $\psi^l$ (resp. twice the order of $\psi^l$) if this order is even (resp. odd). Taking into account that $\sigma \psi^{(p-1)/2}$ commutes with $\rho$, we obtain the groups $\langle \rho,\sigma\psi^{(p-1)/d}\rangle \simeq C_p \rtimes C_d$, for $d$ a divisor of $p-1$, such that $d$ is divisible by 4, and the groups $\langle \rho,\sigma\psi^{2(p-1)/d}\rangle \simeq C_p \rtimes C_d$, for $d$ a divisor of $p-1$, such that $d$ is exactly divisible by 2. Taking into account that $\sigma$ commutes with $\rho\varphi^2$, we obtain the groups $\langle \rho,\sigma\psi^{(p-1)/d}\rangle \simeq C_p \rtimes C_d$, for $d$ an even divisor of $p-1$. We have then $2p$ transitive subgroups of $\Hol(D_{2p})$ isomorphic to the semidirect product $C_{p} \rtimes C_d$, for each even divisor $d$ of $p-1$. This gives case (iv) in the statement.

    We determine now the number of structures by using Corollary \ref{cor}. Taking into account that $|\Aut(N)|=|\Aut(D_{2p})|=p(p-1)$, and $|\Aut(C_{2p})|=|\Aut(C_p)|=p-1$ and that for all groups $G$ described above, an automorphism sending $G'$ to $G'$ must restrict on $G'$ to the identity, we obtain that the number of structures is always 2. $\Box$

\end{enumerate}

\begin{remark}  We note that all integer numbers $g$ with $2\leq g\leq 11$, except $g=8$, are of one of the forms $p, p^2$ or $2p$, with $p$ prime. The prime case has been considered in \cite{Ch1} and \cite{P}. The results obtained by our program have allowed us to intuit the classification of separable extensions of degree $2p$ or $p^2$ with respect to their Hopf Galois character. Since $8$ is of none of these forms, the case $g=8$ is specially interesting and besides it presents a high richness of results.
\end{remark}

\section{Example}
In this section we perform explicitly the bijection from the set of regular subgroups of the symmetric group $S_g$ normalized by $\lambda(G)$ to the set of isomorphism classes of Hopf Galois structures given by Theorem \ref{G-P} for a particular example. As stated above, degree 8 extensions exhibit a high richness of results. We shall examine the case presenting the biggest Hopf algebra isomorphism class.
We consider a Galois extension $L/K$ with Galois group $G=C_2\times C_2\times C_2$. As given in table \ref{d8-ext-Giso}, it has 42 Hopf Galois structures of type $D_{2\cdot 4}$ partitioned in 7 Hopf algebra isomorphism classes of 6 elements each. We will examine in detail one of these classes and determine the corresponding Hopf algebra and Hopf actions. We may write $L=K(\alpha,\beta,\gamma)$, with $\alpha^2, \beta^2, \gamma^2 \in K$ and $G$ is then generated by the automorphisms $a,b,c$ given by

$$\begin{array}{cccr} a:& \alpha & \mapsto & -\alpha \\ & \beta & \mapsto & \beta \\ &\gamma & \mapsto & \gamma \end{array}, \quad
\begin{array}{cccr} b:& \alpha & \mapsto & \alpha \\ & \beta & \mapsto & -\beta \\ &\gamma & \mapsto & \gamma \end{array}, \quad
\begin{array}{cccr} c:& \alpha & \mapsto & \alpha \\ & \beta & \mapsto & \beta \\ &\gamma & \mapsto & -\gamma \end{array}.$$

\noindent The group $C_2\times C_2\times C_2\simeq 8T3$ is given in Magma as the subgroup of the symmetric group generated by $(1,8)(2,3)(4,5)(6,7),(1,3)(2,8)(4,6)(5,7),(1,5)(2,6)(3,7)(4,8)$. If we order the elements in $G$ as $Id, ab, b, ac, c, abc, bc, a$, we have $\lambda(a)=(1,8)(2,3)(4,5)(6,7), \linebreak \lambda(b)=(1,3)(2,8)(4,6)(5,7), \lambda(c)=(1,5)(2,6)(3,7)(4,8)$ and we shall identify $G$ with its image by $\lambda$. The following regular subgroups of $S_8$ are isomorphic to $D_{2\cdot 4}$, normalized by $G$ and mutually $G$-isomorphic.

$$\begin{array}{c} N_1=\langle s_1=(1,8)(2,3)(4,5)(6,7), r_1=(1,6,5,2),(3,4,7,8) \rangle, \\ N_2=\langle s_2=(1,8)(2,3)(4,5)(6,7), r_2=(1,4,7,2),(3,6,5,8) \rangle, \\ N_3=\langle s_3=(1,6)(2,5)(3,4)(7,8), r_3=(1,4,5,8),(2,3,6,7) \rangle, \\ N_4=\langle s_4=(1,2)(3,8)(4,7)(5,6), r_4=(1,4,3,6),(2,5,8,7) \rangle, \\ N_5=\langle s_5=(1,6)(2,5)(3,4)(7,8), r_5=(1,2,3,8),(4,5,6,7) \rangle, \\ N_6=\langle s_6=(1,4)(2,7)(3,6)(5,8), r_6=(1,6,7,8),(2,3,4,5) \rangle. \end{array}$$

We check that $ar_ia=r_i^3, br_ib=r_i, cr_ic=r_i, as_ia=s_i, bs_ib=s_i, cs_ic=s_i, 1\leq i \leq 6$, hence $N_i$ is normalized by $G$, for $1\leq i \leq 6$ and $s_i \mapsto s_j, r_i \mapsto r_j$ defines a $G$-isomorphism from $N_i$ to $N_j, 1\leq i,j \leq 6$.

By computation, we obtain that the Hopf algebra corresponding to $N_i$ is the $K$-Hopf algebra with basis $1, r_i+r_i^3, r_i^2, \alpha (r_i-r_i^3), s_i, s_ir_i+s_ir_i^3, s_ir_i^2, \alpha (s_ir_i-s_ir_i^3)$ and the Hopf actions are given by

$$\begin{array}{lllllll} r_1 \mapsto ab, & r_1^2 \mapsto c, & r_1^3 \mapsto abc, & s_1 \mapsto a, & s_1r_1 \mapsto bc, & s_1r_1^2 \mapsto ac, & s_1r_1^3 \mapsto b \\
r_2 \mapsto ab, & r_2^2 \mapsto bc, & r_2^3 \mapsto ac, & s_2 \mapsto a, & s_2r_2 \mapsto c, & s_2r_2^2 \mapsto abc, & s_2r_2^3 \mapsto b \\
 r_3 \mapsto a, & r_3^2 \mapsto c, & r_3^3 \mapsto ac, & s_3 \mapsto abc, & s_3r_3 \mapsto b, & s_3r_3^2 \mapsto ab, & s_3r_3^3 \mapsto bc \\
 r_4 \mapsto abc, & r_4^2 \mapsto b, & r_4^3 \mapsto ac, & s_4 \mapsto ab, & s_4r_4 \mapsto bc, & s_4r_4^2 \mapsto a, & s_4r_4^3 \mapsto c \\
r_5 \mapsto ab, & r_5^2 \mapsto b, & r_5^3 \mapsto a, & s_5 \mapsto abc, & s_5r_5 \mapsto c, & s_5r_5^2 \mapsto ac, & s_5r_5^3 \mapsto bc \\
 r_6 \mapsto abc, & r_6^2 \mapsto bc, & r_6^3 \mapsto a, & s_6 \mapsto ac, & s_6r_6 \mapsto b, & s_6r_6^2 \mapsto ab, & s_6r_6^3 \mapsto c
\end{array}
$$

A different explicit example can be found in \cite{S}, Example 5.3.1.

\section{Program output}

We present the results obtained for separable field extensions of degree $8$ in Tables \ref{d8-ext}, \ref{d8-ext-cont} and \ref{d8-ext-Giso} in the appendix. We denote by $kTi$ the $i$th transitive group of degree $k$ called by TransitiveGroup$(k,i)$ in the Magma program. In Tables \ref{d8-ext} and \ref{d8-ext-Giso}, for each regular subgroup of $S_8$ (i.e. for $i=1,\dots,5$), we give the name of the abstract group of order 8 isomorphic to it. For the names of some of the remaining groups, the reader may consult Table 8A in \cite{Bu}. In Tables \ref{d8-ext} and \ref{d8-ext-cont}, for each transitive group $G$ of degree $8$ and each group $N$ of order $8$, we give the total number $T$ of Hopf Galois structures of type $N$ for a separable field extension $L/K$ of degree $8$ such that the Galois group of the normal closure $\widetilde{L}$ over $K$ is isomorphic to $G$. Moreover, we give the number a-c of those which are almost classically Galois, the number BC of those for which the Galois correspondence is bijective and the number G-i of Hopf algebra isomorphism classes in which the Hopf Galois structures are partitioned. In particular the difference BC minus a-c gives the number of non almost classically Galois Hopf Galois structures for which the Galois correspondence is bijective. The transitive groups $G$ such that the corresponding field extension $L/K$ has no Hopf Galois structure are not included in the table.

\vspace{0.3cm}
We note that the field extension with smallest degree having a non almost classically Galois Hopf Galois structure with bijective Galois correspondence is a Galois extension of degree 4 with Galois group $C_4$ and Hopf Galois structure of type $C_2\times C_2$. The non-Galois extension with smallest degree having this property is a separable extension of degree 6 whose Galois closure has group $6T5$ and the Hopf Galois structure is of type $S_3$. The field extension with smallest degree having non-isomorphic Hopf Galois structures with isomorphic Hopf algebras is a Galois extension of degree 6 with Galois group the symmetric group $S_3$ for which the three Hopf Galois structures of cyclic type $C_6$ have underlying isomorphic Hopf algebras.

\vspace{0.3cm}
In table \ref{d8-ext-Giso} we give the distribution of Hopf Galois structures in Hopf algebra isomorphism classes for transitive groups of degree 8 having some class with more that one element. For example, in the cell corresponding to $G=N=C_4\times C_2, 10=5\times 1+1 \times 2+1 \times 3$ means that for a Galois extension with Galois group $C_4\times C_2$ there are 10 Hopf Galois structures of type $C_4\times C_2$ which are distributed in 5 classes with 1 element, 1 class with 2 elements and 1 class with 3 elements.

\vspace{0.3cm}
Table \ref{fig} is a compendium of the computation results. In it we give for every degree $g$ the total number of transitive groups of degree $g$ and the number Max of transitive groups of degree $g$ whose order does not exceed the order of the holomorphs of all the groups of order $g$; the number of possible types of Hopf Galois structures; the total number of Hopf Galois structures and the number of the almost classically Galois ones; the number of Hopf Galois structures with bijective Galois correspondence and the number of those which are not almost classically Galois; the number of Hopf algebra isomorphism classes in which the Hopf Galois structures are partitioned (which correspond to $G$-isomorphism classes of the corresponding regular groups $N$) and the number of those for Galois extensions (i.e. when $G'=\Gal(\wL/L)$ is trivial); and finally the execution times in seconds and the memory used in megabytes. We note that the presented program is very efficient up to degree 11. One may observe in particular that the computation for degree 8, which gives a large number of Hopf Galois structures, takes only about 17 seconds. The memory used reaches 160 megabytes for degree 11.

\vspace{0.3cm}
The reader may find in \cite{CS} for each degree $g$ up to 11, the output of the program containing in particular, the precise description of the regular subgroups $N$ of $S_g$ corresponding to the Hopf Galois structures as well as tables summarizing these results.

\section{Conclusions}

The elaboration of the program presented allows to determine all Hopf Galois structures of separable field extensions of a given degree up to degree 11. Such a determination has been obtained by theoretic tools only for prime degree extensions. Besides, a careful inspection of the data provided by the program has led us to obtain several theoretic results. In proposition \ref{psqua} we prove a partial result concerning Hopf Galois structures of separable field extensions of degree $p^2$, for $p$ an odd prime, and describe in theorem \ref{p^2} those of cyclic type. This result came up from the output of the program for degree 9 extensions. In theorems \ref{2p1} and \ref{2p2} we determine Hopf Galois structures of separable field extensions of degree $2p$, for $p$ an odd prime. This result was suggested by the output for degree 6 and 10 extensions. The results obtained by the program have given us the intuition about the general behaviour for the infinite families of extensions discussed in the theorems.  We specially highlight the richness of results obtained in the degree 8 case.

\section*{Acknowledgments} We are grateful to Anna Rio and Montserrat Vela for valuable discussions on the subject of this paper and to Joan Nualart and Pawe\l \ Bogdan for their help with the Magma program.

\newpage
\addtolength{\textwidth}{4cm}
\begin{landscape}
\section*{Appendix - Tables}

\begin{table}[h]
\centering
\caption{Degree 8 extensions}
\label{d8-ext}
\vspace{0.2 cm}
\begin{tabular}{|c||c|c|c|c||c|c|c|c||c|c|c|c||c|c|c|c||c|c|c|c|}
\hline
\multicolumn{1}{|c||}{} & \multicolumn{20}{|c|}{\bf Hopf Galois structures} \\
\hline \hline
 \multicolumn{1}{|c||}{\bf Galois } & \multicolumn{4}{|c||}{Type $C_8$} & \multicolumn{4}{|c||}{Type $C_4\times C_2$} & \multicolumn{4}{|c||}{Type $C_2\times C_2\times C_2$} & \multicolumn{4}{|c||}{Type $D_{2\cdot4}$} & \multicolumn{4}{|c|}{Type $Q_8$}\\ \cline{2-21}  \multicolumn{1}{|c||}{\bf group } & T & a-c &  BC & G-i & T & a-c  & BC & G-i & T & a-c  & BC & G-i & T & a-c  & BC & G-i & T & a-c  & BC & G-i \\
\hline
\multicolumn{1}{|c||}{$8T1\simeq C_8$}                     & 2 &1&2&2& 0 &-&-&-     & 0 &-&-&-               & 2 &0&2&2     & 2 &0&2&2    \\ \hline
\multicolumn{1}{|c||}{$8T2\simeq C_4\times C_2$}           & 4 &0&0&2 & 10  &1&1&7 & 4 &0&1&4                 & 6&0&2&5    & 2  &0&0&2   \\ \hline
\multicolumn{1}{|c||}{$8T3\simeq (C_2)^3$} & 0&-&-&-  & 42 &0&0&28    & 8 &1&1&8              & 42 &0&0&7   & 14 &0&0&7   \\ \hline
\multicolumn{1}{|c||}{$8T4\simeq D_{2\cdot4}$}             & 2 &0&0&1 & 14 &0&0&9     & 6  &0&0&4             & 6 &1&1&4   & 2 &0&0&2    \\ \hline
\multicolumn{1}{|c||}{$8T5\simeq Q_8$}                     & 6&0&0&3  & 6  &0&6&3     & 2 &0&2&1              & 6 &0&6&6   & 2 &1&2&2  \\ \hline
\multicolumn{1}{|c||}{$8T6$}                    & 2 &1&2&2   & 0  &-&-&-             & 0 &-&-&-               & 2 &1&2&2   & 2 &0&2&2    \\ \hline
\multicolumn{1}{|c||}{$8T7$}                     & 2&2&2&2     & 0 &-&-&-              & 0 &-&-&-             & 2&0&2&2    & 2&0&2&2     \\ \hline
\multicolumn{1}{|c||}{$8T8$}                     & 2 &1&2&2    & 0 &-&-&-              & 0 &-&-&-             & 2&0&2&2    & 2 &1&2&2    \\ \hline
\multicolumn{1}{|c||}{$8T9$}                     & 0 &-&-&-    & 10&1&1&9           & 4&1&1&4                 & 6 &2&2&5   & 2 &0&0&2    \\ \hline
\multicolumn{1}{|c||}{$8T10$}                     & 0 &-&-&-    & 6 &2&3&6          & 4 &0&1&4       & 0 &-&-&-            & 0 &-&-&-    \\ \hline
\multicolumn{1}{|c||}{$8T11$}                     & 2 &0&0&1   & 6  &2&6&5          & 2 &0&2&2       & 6  &1&6&6           & 2 &1&2&2  \\ \hline
\multicolumn{1}{|c||}{$8T12$}                     & 0&-&-&-     & 0&-&-&-               & 2 &0&2&1   & 0 &-&-&-           & 2  &1&2&2   \\ \hline
\multicolumn{1}{|c||}{$8T13$}                     & 0 &-&-&-    & 0  &-&-&-             & 2 &1&1&2  & 0 &-&-&-           & 2 &0&0&1    \\ \hline
\multicolumn{1}{|c||}{$8T14$}                     & 0 &-&-&-   & 0&-&-&-              & 4 &0&1&3       & 0&-&-&-            & 0 &-&-&-    \\ \hline
\multicolumn{1}{|c||}{$8T15$}                     & 2 &2&2&2    & 0&-&-&-            & 0&-&-&-         & 2 &1&2&2         & 2&1&2&2     \\ \hline
\multicolumn{1}{|c||}{$8T16$}                     & 0&-&-&-    & 0 &-&-&-              & 0 &-&-&-      & 2 &0&2&2         & 2&0&2&2     \\ \hline
\multicolumn{1}{|c||}{$8T17$}                     & 0&-&-&-     & 0&-&-&-               & 0 &-&-&-     & 2  &1&2&2           & 2&1&2&2     \\ \hline
\multicolumn{1}{|c||}{$8T18$}                     & 0&-&-&-     & 6 &3&3&6              & 4 &1&1&4     & 0 &-&-&-            & 0 &-&-&-    \\ \hline

\end{tabular}
\end{table}

\begin{table}[t]
\centering
\caption{Degree 8 extensions (cont.)}
\label{d8-ext-cont}
\vspace{0.2 cm}
\begin{tabular}{|c||c|c|c|c||c|c|c|c||c|c|c|c||c|c|c|c||c|c|c|c|}
\hline
\multicolumn{1}{|c||}{} & \multicolumn{20}{|c|}{\bf Hopf Galois structures} \\
\hline \hline
 \multicolumn{1}{|c||}{\bf \quad \, Galois \quad \, \, } & \multicolumn{4}{|c||}{Type $C_8$} & \multicolumn{4}{|c||}{Type $C_4\times C_2$} & \multicolumn{4}{|c||}{Type $C_2\times C_2\times C_2$} & \multicolumn{4}{|c||}{Type $D_{2\cdot4}$} & \multicolumn{4}{|c|}{Type $Q_8$}\\ \cline{2-21}  \multicolumn{1}{|c||}{\bf group } & T & a-c &  BC & G-i & T & a-c  & BC & G-i & T & a-c  & BC & G-i & T & a-c  & BC & G-i & T & a-c  & BC & G-i \\
\hline
\multicolumn{1}{|c||}{ $8T19$}                     & 0 &-&-&-    & 2 &1&2&2              & 2 &1&2&2     & 0&-&-&-            & 0 &-&-&-    \\ \hline
\multicolumn{1}{|c||}{$8T20$}                     & 0&-&-&-     & 2 &0&2&2              & 2 &0&2&2     & 0 &-&-&-            & 0 &-&-&-    \\ \hline
\multicolumn{1}{|c||}{$8T22$}                     & 0&-&-&-     & 6 &6&6&6             &2 &2&2&2       & 6  &6&6&6            & 2  &2&2&2  \\ \hline
\multicolumn{1}{|c||}{$8T23$}                     & 0&-&-&-     & 0 &-&-&-              & 0 &-&-&-     & 0 &-&-&-          & 2 &1&2&2     \\ \hline
\multicolumn{1}{|c||}{$8T24$}                     & 0&-&-&-    & 0&-&-&-               & 2 &1&1&2     & 0 &-&-&-           & 0 &-&-&-    \\ \hline
\multicolumn{1}{|c||}{$8T25$}                     & 0  &-&-&-   & 0  &-&-&-             & 1  &1&1&1    & 0 &-&-&-           & 0 &-&-&-   \\ \hline
\multicolumn{1}{|c||}{$8T26$}                     & 0 &-&-&-    & 0 &-&-&-              & 0 &-&-&-     & 2  &2&2&2           & 2  &2&2&2 \\ \hline
\multicolumn{1}{|c||}{$8T29$}                     & 0 &-&-&-    & 2 &2&2&2              & 2  &2&2&2      & 0 &-&-&-           & 0 &-&-&-    \\ \hline
\multicolumn{1}{|c||}{$8T32$}                     & 0  &-&-&-   & 0  &-&-&-            & 2  &2&2&2      & 0 &-&-&-           & 2  &2&2&2   \\ \hline
\multicolumn{1}{|c||}{$8T33$}                     & 0 &-&-&-    & 0 &-&-&-              & 1 &1&1&1       & 0 &-&-&-           & 0 &-&-&-   \\ \hline
\multicolumn{1}{|c||}{$8T34$}                     & 0 &-&-&-    & 0 &-&-&-              & 3 &0&3&3       & 0 &-&-&-           & 0 &-&-&-    \\ \hline
\multicolumn{1}{|c||}{$8T36$}                     & 0 &-&-&-    & 0 &-&-&-             & 1 &1&1&1       & 0 &-&-&-          & 0 &-&-&-   \\ \hline
\multicolumn{1}{|c||}{$8T37$}                     & 0&-&-&-     & 0 &-&-&-              & 2 &0&2&2      & 0 &-&-&-           & 0&-&-&-     \\ \hline
\multicolumn{1}{|c||}{$8T39$}                     & 0 &-&-&-    & 0 &-&-&-              & 2&2&2&2                         & 0  &-&-&-          & 0&-&-&-     \\ \hline
\multicolumn{1}{|c||}{$8T40$}                     & 0 &-&-&-    & 0 &-&-&-              & 0 &-&-&-        &0 &-&-&-            & 2&2&2&2     \\ \hline
\multicolumn{1}{|c||}{$8T41$}                     & 0 &-&-&-    & 0  &-&-&-             & 1 &1&1&1                        & 0 &-&-&-           & 0 &-&-&-    \\ \hline
\multicolumn{1}{|c||}{$8T48$}                     & 0 &-&-&-    & 0  &-&-&-             & 1 &1&1&1      & 0 &-&-&-           & 0 &-&-&-    \\ \hline
\end{tabular}
\end{table}

\clearpage

\addtolength{\hoffset}{-2cm}
\addtolength{\textwidth}{4cm}

\begin{table}[t]
\centering
\caption{Hopf algebra isomorphism classes for degree 8 extensions}
\label{d8-ext-Giso}
\vspace{0.2 cm}
\begin{tabular}{|c||c|c|c|c|c|}
\hline
\multicolumn{1}{|c||}{} & \multicolumn{5}{|c|}{\bf Isomorphism classes} \\
\hline \hline
 \multicolumn{1}{|c||}{\bf{Galois group}} & \multicolumn{1}{|c|}{Type $C_8$} & \multicolumn{1}{|c|}{Type $C_4\times C_2$} & \multicolumn{1}{|c|}{Type $C_2\times C_2\times C_2$} & \multicolumn{1}{|c|}{Type $D_{2\cdot4}$} & \multicolumn{1}{|c|}{Type $Q_8$}\\ \hline
\multicolumn{1}{|c||}{$8T2\simeq C_4\times C_2$}           & $4=2\times 2$ & $10=5\times 1+1 \times 2+1 \times 3$ & $4=4\times 1$&$6=4\times 1+1 \times2$& $2=2\times 1$  \\ \hline
\multicolumn{1}{|c||}{$8T3\simeq (C_2)^3$} & 0  & $42=21\times 1+7\times 3$    & $8=8\times 1$ &$42=7\times 6$ & $14=7\times 2$   \\ \hline
\multicolumn{1}{|c||}{$8T4\simeq D_{2\cdot4}$}  & $2=1\times 2$ & $14=4\times 1+5\times 2$ & $6=2\times 1+2\times 2$ & $6=3\times 1+1\times 3$ & $2=2\times 1$   \\ \hline
\multicolumn{1}{|c||}{$8T5\simeq Q_8$} &$6=3\times 2$ &$6=3\times 2$ & $2=1\times 2$ &$6=6\times 1$ &$2=2\times 1$ \\ \hline
\multicolumn{1}{|c||}{$8T9$} &$0$&$10=8\times 1+1\times 2$&$4=4\times 1$&$6=4\times 1+1\times 2$&$2=2\times 1$    \\ \hline
\multicolumn{1}{|c||}{$8T11$}  &$2=1\times 2$&$6=4\times 1+1\times 2$&$2=2\times 1$&$6=6\times 1$&$2=2\times 1$    \\ \hline
\multicolumn{1}{|c||}{$8T12$} &$0$&$0$&$2=1\times 2$&$0$&$2=2\times 1$    \\ \hline
\multicolumn{1}{|c||}{$8T13$}   & $0$ &$0$&$2=2\times 1$&$0$&$2=1\times 2$    \\ \hline
\multicolumn{1}{|c||}{$8T14$}                     & $0$ &$0$&$4=2\times 1+1\times 2$&$0$    & $0$    \\ \hline
\end{tabular}
\end{table}

\begin{table}[ht]
\centering
\caption{Summary of results}
\label{fig}
\vspace{0.2 cm}
\begin{tabular}{|c||c|c||c||c|c||c|c||c|c||c||c|}
\hline
\multicolumn{1}{|c||}{\bf Degree} & \multicolumn{2}{|c||}{\bf Transitive Groups} &\multicolumn{1}{|c||}{\bf Types} &\multicolumn{2}{|c||}{\bf HG struct.} &\multicolumn{2}{|c||}{\bf BC} &\multicolumn{2}{|c||}{\bf $G$-iso} & \multicolumn{1}{|c||}{\bf Execution time} &\multicolumn{1}{|c|}{\bf Memory used} \\
\cline{2-3} \cline{5-10}
 \multicolumn{1}{|c||}{} &  \multicolumn{1}{|c|}{\quad Total\quad} & \multicolumn{1}{|c||}{Max} & \multicolumn{1}{|c||}{} & \multicolumn{1}{|c|}{Total} & \multicolumn{1}{|c||}{a-c}& \multicolumn{1}{|c|}{Total} & \multicolumn{1}{|c||}{not a-c}& \multicolumn{1}{|c|}{Total}& \multicolumn{1}{|c||}{Galois}&\multicolumn{1}{|c||}{(s)}&\multicolumn{1}{|c|}{(MB)}\\ \hline \hline
 2&1&1&1&1&1&1&0&1&1&$\approx 1$ & $\approx 10$  \\ \hline
3&2&2&1&2&2&2&0&2&1&$\approx 1$ & $\approx 10$  \\ \hline
4&5&5&2&10&6&7&1&10&6&$\approx 1$ & $\approx 11$ \\ \hline
5&5&3&1& 3&3&3&0&3&1&$\approx 1$ &$\approx 11$ \\ \hline
6&16&10&2&15&7&9&2&13&6&$\approx 2$ & $\approx 11$\\ \hline
7& 7&4&1&4&4&4&0&4&1&$\approx 1$ &$\approx 11$ \\ \hline
8&50&48& 5&348&74&147&73&262&111&$\approx 17$&$\approx 40$ \\ \hline
9& 34&26&2&38&26&28&2&33&8&$\approx 10$& $\approx 16$ \\ \hline
10&45&21&2&27&11&17&6&23&6&$\approx 160$&$\approx 45$ \\ \hline
11& 8&4&1&4&4&4&0&4&1&$\approx 90$ & $\approx 160$ \\ \hline
\end{tabular}
\end{table}

\end{landscape}
\end{document}